\documentclass[12pt]{amsart}
\usepackage{amsfonts}
\usepackage{amssymb}
\theoremstyle{definition}

\theoremstyle{remark}

\newcommand{\C}{\mathbb C}
\newcommand{\D}{\mathbb D}
\newcommand{\Sph}{\mathbb S}

\newcommand{\X}{\mathfrak X}

\newcommand{\F}{\mathfrak F}
\newcommand{\ds}{\displaystyle}

\begin{document}

\centerline{\large\bf ALMOST \ K\"AHLER \ MANIFOLDS \ OF \ CONSTANT}
\centerline{\large\bf ANTIHOLOMORPHIC \ SECTIONAL \ CURVATURE
\footnote{\it SERDICA Bulgaricae mathematicae publicationes. Vol. 9, 1983, p. 372-376.}}

\vspace{0.2in}
\centerline{\large OGNIAN T. KASSABOV}

\vspace{0.3in}
{\sl We prove that an $AK_2$-manifold of dimension $2m\ge6$ and of pointwise
constant anti\-holomorphic sectional curvature is either a 6-dimensional manifold
of constant negative sectional curvature or is locally isometric to $\C^n$, 
$\C\mathbb P^n$ or $\C\D^n$.}

\vspace{0.1in}
{\bf 1. Introduction.} Let $M$ be an almost Hermitian manifolds with
metric tensor $g$, almost complex structure $J$ and covariant
differentiation $\nabla$. 

If $\nabla J=0$, or $(\nabla_XJ)X=0$ or \vspace{-.2cm}
$$
	g((\nabla_XJ)Y,Z)+g((\nabla_YJ)Z,X)+g((\nabla_ZJ)X,Y)=0 \,, \leqno (1.1)
$$
\vspace{-.2cm}
or 
\vspace{-.2cm}
$$
	(\nabla_XJ)Y+(\nabla_{JX}J)JY=0    \leqno (1.2)
$$
for all $X,\,Y,\,Z \in \X(M)$, $M$ is said to be a K\"ahler, or nearly
K\"ahler, or almost K\"ahler, or quasi K\"ahler manifold, respectively. The
corresponding classes of manifolds are denoted by $K$, $NK$, $AK$ and $QK$,
respectively. 

One consider the following identities:

1) $R(X,Y,Z,U)=R(X,Y,JZ,JU)$,

2) $R(X,Y,Z,U)=R(X,Y,JZ,JU)+R(X,JY,Z,JU)+R(JX,Y,Z,JU)$,

3) $R(X,Y,Z,U)=R(JX,JY,JZ,JU)$,

\noindent where $R$ is the curvature tensor for $M$. If $M$ has the identity $i$,
it is said to belong to the class $AH_i$. Then $AH_1 \subset AH_2 \subset AH_3$.
For a given class $L$ of almost Hermitian manifolds one denotes $L_i = L \cap AH_i$.
It is known, that $K=K_1$, $NK=NK_2$, $NK\cap AK=K$. The inclusions 
$K\subset NK\subset QK_2$, $K\subset AK_2 \subset QK_2$ are strict [5].

A 2-plane $\alpha$ in a tangent space $T_p(M)$ is said to be holomorphic 
(respectively, antiholo\-morphic) if $\alpha = J\alpha$ (respectively, 
$\alpha \perp J\alpha$). If for each point $p\in M$ the curvature of an 
arbitrary holomorphic (respectively, antiholomorphic) 2-plane $\alpha$ in
$T_p(M)$ doesn't depend on $\alpha$, $M$ is said to be of pointwise constant
holomorphic (respectively, antiholomorphic) sectional curvature.

Suppose $M \in NK$ and dim\,$M=2m \ge 6$. If $M$ has pointwise constant
holomorphic sectional curvature it is locally isometric to $\C^n$, 
$\C\mathbb P^n$, $\C\D^n$ or $\Sph^6$ [4]. The corresponding result for the 
antiholomorphic case in proved in [2].

In section 3 we prove an analogous theorem for $AK_2$-manifolds of pointwise
constant antiholomorphic sectional curvature.

\vspace{0.1in}
{\bf 2. Preliminaries.} Let $M$ be a $2m$-dimensional almost Hermitian manifold. 
We denote by $R$, $S$ and $\tau(R)$ the curvature tensor, the Ricci tensor and 
the scalar curvature of $M$, respectively. One defines also a tensor $S'$ and 
a function $\tau'(R)$ by
$$
	S'(X,Y)=\sum_{i=1}^{2m} R(X,E_i,JE_i,JY)  \,,
	\qquad \tau'(R)=\sum_{i=1}^{2m} S'(E_i,E_i) \,,
$$
where $\{ E_i;\, i=1,\hdots,2m \}$ is a local orthonormal frame field. From the second 
Bianchi identity one obtains
$$
	\sum_{i=1}^{2m} (\nabla_{E_i}R)(X,Y,Z,E_i)=(\nabla_XS)(Y,Z)-(\nabla_YS)(X,Z) \,, \leqno (2.1)
$$
$$
	\sum_{i=1}^{2m} (\nabla_{E_i}S)(X,E_i)=\frac12 X(\tau(R)) \,. \leqno (2.2)
$$

For $M \in AK_2$ the following identities hold (see e.g. [1], [5]):
$$
	R(X,Y,Z,U)-R(X,Y,JZ,JU)=\frac12 
	        g((\nabla_XJ)Y- (\nabla_YJ)X,(\nabla_ZJ)U-(\nabla_UJ)Z) \,,   \leqno (2.3)
$$
$$
	2(\nabla_X(S-S'))(Y,Z)=(S-S')((\nabla_XJ)Y,JZ)+(S-S')(JY,(\nabla_XJ)Z) \,. \leqno (2.4)
$$

Assume $M\in AH_3$ and dim\,$M=2m$. If $M$ is of pointwise constant antiholomorphic
sectional curvature $\nu$, then [2]:
$$
	R=\frac 16 \psi(S) +\nu\pi_1-\frac{2m-1}3\nu\pi_2 \,,  \leqno (2.5)
$$
$$
	(m+1)S-3S'=\frac1{2m} \{ (m+1)\tau(R)-3\tau'(R) \} g \,,   \leqno (2.6)
$$
where
$$
	\pi_1(x,y,z,u)=g(x,u)g(y,z)-g(x,z)g(y,u)  \,,
$$
$$
	\begin{array}{r}
		\psi(Q)(x,y,z,u)=g(x,Ju)Q(y,Jz)-g(x,Jz)Q(y,Ju)-2g(x,Jy)Q(z,Ju)   \\
		+g(y,Jz)Q(x,Ju)-g(y,Ju)Q(x,Jz)-2g(z,Ju)Q(x,Jy)
	\end{array}
$$
for a tensor $Q$ of type (0,2) and $\ds\pi_2=\frac12\psi(g)$. According to (2.5)
$M \in AH_2$. In the case $m>2$, $\nu$ is a global constant [6]. If moreover 
$M \in QK_3$, $\tau(R)$ and $\tau'(R)$ are also global constants [3].

Consequently, if $ M \in AK_3$, dim\,$M=2m\ge 6$ and $M$ is of pointwise
constant antiholo\-morphic sectional curvature, from (2.2), (2.4) and (2.6)
we obtain:
$$
	\sum_{i=1}^{2m} (\nabla_{E_i}S)(X,E_i)=0 \,, \leqno (2.7)
$$
$$
	2(\nabla_XS)(Y,Z)=S((\nabla_XJ)Y,JZ)+S(JY,(\nabla_XJ)Z) \,. \leqno (2.8)
$$

In section 3 we shall use the following lemma:

{\bf Lemma.} {\it Let $M\in AK_3$ with {\rm dim}$M=2m \ge 6$. If the curvature
tensor of $M$ has the form
$$
	R=f\pi_1 + h \pi_2 \,,    \leqno (2.9)
$$
where $ f,\,h \in \F(M)$, then either $M$ is a 6-dimensional manifold of
constant negative sectional curvature or $M$ is a K\"ahler manifold of constant
holomorphic sectional curvature.}

{\bf Proof.} If $M$ is not a K\"ahler manifold of constant holomorphic sectional 
curvature, it is of constant sectional curvature [8]. If $m>3$ this is impossible [7].
On the other hand, as easily follows from (2.3), if a non K\"ahler almost
K\"ahler manifold is of constant curvature $f$, then $f<0$.

\vspace{0.1in}
{\bf 3. The classification theorem for the class $AK_3$.}

{\bf Theorem.} {\it Let $M \in AK_3$ and {\rm dim}\,$M=2m \ge 6$. If $M$ is of poitwise 
constant antiho\-lomortphic sectional curvature $\nu$, then either $M$ is a
6-dimensional manifold of constant sectional curvature $\nu<0$, or $M$ is locally
isometric to one of the following manifolds:

a) the complex Euclidian space $\C^n$, 

b) the complex projective space $\C\mathbb P^n(4\nu)$;

c) the complex hyperbolic space $\C\D^n(4\nu)$. }

{\bf Proof.} First we prove that the Ricci tensor is parallel. Let $ p \in M$, 
$ x,\,y \in T_p(M)$. According to the second Bianchi identity
$$
	(\nabla_xR)(Jx,y,y,Jx)+(\nabla_{Jx}R)(y,x,y,Jx)+(\nabla_yR)(x,Jx,y,Jx)=0 \,.   \leqno (3.1)
$$

Let $ \{ e_i,\,Je_i; \, i=1,\hdots,m \}$ be an orthonormal basis of $T_p(M)$
such that $S(e_i)=\lambda_ie_i$ (and hence $S(Je_i)=\lambda_iJe_i$), $i=1,\hdots,m$.
In (3.1) we put $x=e_i$, $y=e_j$ ($i\ne j$) and using (1.2), (2.5) and $\nu=$\,const. 
we obtain
$$
	\{ 3\lambda_i+\lambda_j-4(2m-1)\nu \} g((\nabla_{e_j}J)e_j,Je_i)=0 \,.  \leqno(3.2)
$$
Analogously, from 
$$
	(\nabla_{Je_i}R)(e_k,e_j,e_j,Je_k)+(\nabla_{e_k}R)(e_j,Je_i,e_j,Je_k)+(\nabla_{e_j}R)(Je_i,e_k,e_j,Je_k)=0 
$$
for $i\ne j\ne k \ne i $ we find
$$
	\{ 2\lambda_k+\lambda_i+\lambda_j-4(2m-1)\nu \} g((\nabla_{e_j}J)e_j,Je_i)=0 \,.  \leqno(3.3)
$$
If for some $s$ $(\nabla_{e_s}J)e_s \ne 0$, without loss of generality we assume that 
$g((\nabla_{e_s}J)e_s,Je_i) \ne 0$ for some $i\ne s$ and from (3.2) and (3.3) if follows:
$$
	\begin{array}{l}\vspace{0.2cm}
		3\lambda_i+\lambda_s-4(2m-1)\nu = 0 \,, \\
		2\lambda_k+\lambda_i+\lambda_s-4(2m-1)\nu= 0 \,.
	\end{array}   \leqno (3.4)
$$
Hence $ \lambda_i = \lambda_k $ for $i, k=1,\hdots,m$; $i,\,k \ne s$. If there exists 
$ i \ne s$ such that  $(\nabla_{e_i}J)e_i \ne 0$ in the same way we conclude that 
$\lambda_k=\lambda_s$ for $ k\ne i,s$ and consequently $\lambda_i=\lambda_j$, $i,\,j=1,\hdots,m$.
Hence, using (2.8), we find $(\nabla_xS)(y,z)=0$ for all $x,\,y,\,z\in T_p(M)$. Let 
$$
	(\nabla_{e_i}J)e_i=(\nabla_{Je_i}J)Je_i=0    \leqno (3.5)
$$
for $i=1,\hdots,m$, $i\ne s$. Using (2.1), (2.5), (2.7), (1.1), (2.8) and 
$\nu=$\,const. we obtain
$$
	\begin{array}{r}\vspace{0.2cm}
		4S(Jx,(\nabla_zJ)y) - 4S(Jy,(\nabla_zJ)x) + 5 S(Jx,(\nabla_yJ)z) - 5 S(Jy,(\nabla_xJ)z) \ \  \\
		-S((\nabla_xJ)y,Jz)+S((\nabla_yJ)x,Jz)- 12(2m-1) \nu g(Jx,(\nabla_zJ)y)=0  \,.
	\end{array}   \leqno (3.6)
$$
In (3.6) we put $x=e_i$, $y=e_j$, $z=e_k$:
$$
	\begin{array}{r}\vspace{0.2cm}
		\{ 4\lambda_i+4\lambda_j-\lambda_k-12(2m-1)\nu \} g(Je_i,(\nabla_{e_k}J)e_j) \ \  \\
		  +5\lambda_i g(Je_i,(\nabla_{e_j}J)e_k)-5\lambda_j g(Je_j,(\nabla_{e_i}J)e_k)=0  
	\end{array}   \leqno (3.7)
$$
for arbitrary $i,\,j,\,k$. Let $i,\,j \ne s$, $k=s$. Then from (3.7), (1.1) and $\lambda_i=\lambda_j$
it follows $\{ 13\lambda_i-\lambda_s-12(2m-1)\nu \} g(Je_i,(\nabla_{e_s}J)e_j)=0$. If
$ g(Je_i,(\nabla_{e_s}J)e_j)\ne 0$ for some $l,\,j\ne s$, then
$$
	13\lambda_i-\lambda_s-12(2m-1)\nu =0
$$
and using (3.4) we find $\lambda_i=\lambda_s$. Hence, $\nabla S=0$ in $p$. Let
$$
	g(e_i,(\nabla_{e_s}J)e_j)=g(Je_i,(\nabla_{e_s}J)e_j)=0     \leqno (3.8)
$$
for $i,\,j\ne s$. Let in (3.7) $i,\,k\ne s$, $j=s$. Using (3.8) and $\lambda_i=\lambda_k$,
we obtain $\{3\lambda_i+9\lambda_s-12(2m-1)\nu\}g(Je_i,(\nabla_{e_k}J)e_s) =0$. If
$g(Je_i,(\nabla_{e_k}J)e_s) \ne 0$ for some $i,\,k \ne s$, we have $\lambda_i=\lambda_s$ and
hence $\nabla S=0$ in $p$. Let 
$$
	g(e_i,(\nabla_{e_k}J)e_s)=g(Je_i,(\nabla_{e_k}J)e_s)=0  \,.     \leqno (3.9)
$$ 
If $m>3$, let $i,\,j,\,k $ be different and $i,\,j,\,k \ne s$. From (3.7) it follows 
$\{\lambda_i-(2m-1)\nu\}g(Je_i,(\nabla_{e_k}J)e_j) =0$. If $g(Je_i,(\nabla_{e_k}J)e_j) \ne 0$ then
$\lambda_i=(2m-1)\nu$ and using (3.4) we obtain $\lambda_i=\lambda_s$. Hence and because
of (3.8), (3.9) we assume that 
$$
	g(e_i,(\nabla_{e_j}J)e_k)=g(Je_i,(\nabla_{e_j}J)e_k)=0
$$
for all different $i,\,j,\,k$. According to these equalities and (3.5)
$$
	(\nabla_{e_i}J)e_j=0  \leqno (3.10)
$$
for $i\ne s$ and for arbitrary $j=1,\hdots,m$. From (2.3) and (3.10) we derive
$R(e_i,e_j,e_j,e_i)-R(e_i,e_j,Je_j,Je_i)=0$ for $i\ne j$; $i,j\ne s$. Hence, 
because of (2.5) and $\lambda_i = \lambda_j$ it follows $\lambda_i=2(m+1)\nu$
for $i \ne s$ and using (3.4) we find $\lambda_s=2(m-5)\nu$. On the other hand,
from (2.3) and (3.10) we obtain
$$
	2R(e_i,e_s,e_s,e_i)-2R(e_i,e_s,Je_s,Je_i)=-g((\nabla_{e_s}J)e_i),(\nabla_{e_s}J)e_i) \,.
$$ 
From the last three equalities and (2.5) we derive
$$
	g((\nabla_{e_s}J)e_i),(\nabla_{e_s}J)e_i)=-4\nu \,.  \leqno (3.11)
$$
On the other hand, because of (3.8), $(\nabla_{e_s}J)e_i=\alpha_ie_s+\beta_iJe_s$,
where $\alpha_i,\,\beta_i$ are real constants, $i=1,\hdots,m$, $i\ne s$. Using again
(2.3), (2.5) and (3.10), we derive
$$
	g((\nabla_{e_s}J)e_i),(\nabla_{e_s}J)e_j)=0\,, \qquad
	g((\nabla_{e_s}J)e_i),(\nabla_{e_s}J)Je_j)=0
$$ 
for $i\ne j$; $i,\,j\ne s$. Consequently, $\alpha_i\alpha_j+\beta_i\beta_j=0$,
$-\beta_i\alpha_j+\alpha_i\beta_j=0$. Hence, by using (3.11) it is easy to find
$\nu = 0$ and $(\nabla_{e_s}J)e_i=0$ for $i=1,\hdots,m$; $i\ne s$. Hence,
$(\nabla_{e_s}J)e_s=0$ which is a contradiction.

Consequently, from $(\nabla_{e_s}J)e_s \ne 0$ for some $s$ it follows $\nabla S=0$
in $p$.

Let $(\nabla_{e_i}J)e_i=0$ for every $i=1,\hdots,m$. Then
$$
	S((\nabla_{e_i}J)e_i,y)=S(e_i,(\nabla_{e_i}J)y)=0   \leqno (3.12)
$$
for every $i=1,\hdots,m$ and for arbitrary $y \in T_p(M)$. Now from (3.1) with
$x=e_i$ by using (2.5), (2.8) and (3.12) we derive
$$
	3S(Je_i,(\nabla_yJ)y)-S(Jy,(\nabla_yJ)e_i)-4(2m-1)\nu g(Je_i,(\nabla_yJ)y)=0 \,.
$$
Hence,
$$
	3S(Jx,(\nabla_yJ)y)-S(Jy,(\nabla_yJ)x)-4(2m-1)\nu g(Jx,(\nabla_yJ)y)=0 
$$
for all $x,\,y\in T_p(M)$. On the other hand, from (3.6) it follows
$$
	3S(Jx,(\nabla_yJ)y)-S(Jy,(\nabla_yJ)x)-2S(Jy,(\nabla_xJ)y)-4(2m-1)\nu g(Jx,(\nabla_yJ)y)=0 \,.
$$
By using (2.8) from the last two identities we find $(\nabla_xS)(y,y)=0$, i.e.
$\nabla S=0$ in $p$.

Consequently, the Ricci tensor is parallel. If $M$ is irreducible it is an Einsteinian 
manifold. Hence, the curvature tensor has the form (2.9) and the assertion follows from 
the Lemma and the classification theorem for K\"ahler manifolds of constant holomorphic
sectional curvature.

Let $M$ be reducible. Then $M$ is locally isometric to a product
$M_1(\mu_1)\times\hdots \times M_k(\mu_k)$, where $S=\mu_ig$ on $M_i(\mu_i)$ and
$\mu_i\ne \mu_j$ for $i\ne j$. Then it is not difficult to prove that $M_i(\mu_i)$
is an $AK_2$-manifold for $i=1,\hdots,m$. If $k=1$, $M$ is an Einsteinian manifold
and the theorem follows. Let $k>1$. By using (2.5) it is easy to see that
$\nu=0$ and $\mu_i=-\mu_j$ for $i \ne j$. If $k>2$ it follows $\mu_i=0$ for
$i=1,\hdots,k$ which is a contradiction. Let $k=2$ and dim\,$M_1(\mu_1) \ge 4$. 
Analogously to (3.2) $\mu_1g(x,(\nabla_yJ)y)=0$ holds good for all 
$ x,\,y \in T_p(M_1(\mu_1))$, $p \in M_1(\mu_1)$. Since $\mu_1=0$ is a contradiction,
we assume $\mu_1 \ne 0$. Then $M_1(\mu_1)$ is a nearly K\"ahler manifold and hence 
it is a K\"ahler manifold of constant holomorphic sectional curvature and of zero 
antiholomorphic sectional curvature. Hence, $\mu_1=0$ which is again a contradiction.

\vspace{0.6in}
\centerline{\large R E F E R E N C E S}

\vspace{0.2in}
\noindent
1. M.  B\,a\,r\,r\,o\,s. Classes de Chern de las $NK$-variedades. Geometria de las 
$AK_2$-varieda-

des. Tesis doctorales. Universidad de Granada, 1977.

\noindent
2. G. G\,a\,n\,\v{c}\,e\,v,  O. K\,a\,s\,s\,a\,b\,o\,v. Nearly K\"ahler manifolds of 
constant antiholomorphic 

sectional curvature. {\it C. R. Acad. bulg. Sci.}, {\bf 35}, 1982, 145-147.

\noindent
3. G. G\,a\,n\,\v{c}\,e\,v,  O. K\,a\,s\,s\,a\,b\,o\,v. Schur's theorem of antiholomorphic type
for quasi K\"ah-

ler manifolds. {\it C. R. Acad. bulg. Sci.}, {\bf 35}, 1982, 307-309.

\noindent
4. A. G\,r\,a\,y.  Classification des vari\'et\'es approximativement k\"ahleriennes de courbure sec-

tionnelle holomorphe constante. {\it C. R. Acad. Sci. Paris, S\'er. A}, {\bf 279}, 1974, 797-800.

\noindent
5. A. G\,r\,a\,y. Curvature identities for Hermitian and almost hermitian manifolds. 
{\it T\^ohoku 

Math. J.,} {\bf 28}, 1976, 601-612.

\noindent
6. O. K\,a\,s\,s\,a\,b\,o\,v. Sur le th\'eor\`eme de F. Schur pour une vari\'et\'e presque 
hermitienne. 

{\it C. R. Acad. bulg. Sci.,} {\bf 35}, 1982, 905-908.

\noindent
7. Z. O\,l\,s\,z\,a\,k. A note on almost K\"ahler manifolds. {\it Bull. Acad. Polon. Sci.,
S\'er. Sci. 

Math. Astr. Phis.}, {\bf 26}, 1978, 139-141.

\noindent
8. F. T\,r\,i\,c\,e\,r\,r\,i, L. V\,a\,n\,h\,e\,c\,k\,e. Curvature tensors on almost Hermitian manifolds.

\ {\it Trans. Amer. Math. Soc.}, {\bf 267}, 1981, 365-398.

\vspace {0.5cm}
\noindent
{\it Center for mathematics and mechanics \ \ \ \ \ \ \ \ \ \ \ \ \ \ \ \ \ \ \ \ \ \ \ \ \ \ \ \ \ \ \ \ \ \
Received 15.12.1981

\noindent
1090 Sofia   \ \ \ \ \ \ \ \ \ \ \ \ \ \ \ \ \  P. O. Box 373}

\end{document}